# Approximation of Fractional Derivatives Via Gauss Integration

S.M. Hashemiparast[*], H. Fallahgoul[†]

*Department of Mathematics, Faculty of Science,
K. N. Toosi University of Technology,
P. O. Box $16765 - 165$, Tehran, Iran*

**Abstract**

In this paper approximations of three classes of fractional derivatives (FD) using modified Gauss integration (MGI) and Gauss-Laguerre integration (GLI) are considered. The main solutions of these fractional derivatives depend on inverse of Laplace transforms, which are handled by these procedures. In the modified form of integration the weights and nodes are obtained by means of a difference equation, which gives a proper approximation form for the inverse of Laplace transform and hence the fractional derivatives. Theorems are established to indicate the degree of exactness and boundary of the error of the solutions. Numerical examples are given to illuminate the results of the application of these methods.

[*]hashemiparast@kntu.ac.ir
[†]hfallahgoul@dena.kntu.ac.ir